\documentclass[10pt,a4paper]{article}

\usepackage[latin1]{inputenc}
\usepackage[T1]{fontenc}
\usepackage[english]{babel}
\usepackage{amsfonts}
\usepackage{euscript,amsthm}

\newtheorem{conj}{Conjecture}[section]
\newtheorem{theo}[conj]{Theorem}
\newtheorem{prop}{Proposition}[section]
\newtheorem{coro}[prop]{Corollary}
\newtheorem{lemm}{Lemma}[section]

\newtheorem{defi}{Definition}[section]
\newtheorem{nota}[defi]{Notation}

\begin{document}



\newcommand{\N}{\mathbb{N}}
\newcommand{\Z}{\mathbb Z}
\newcommand{\R}{\mathbb{R}}
\newcommand{\Q}{\mathbb{Q}}
\newcommand{\C}{\mathbb{C}}
\renewcommand{\H}{\mathbb{H}}
\renewcommand{\o}{\mathbb{O}}

\renewcommand{\a}{{\cal A}}
\newcommand{\az}{\a_\Z}
\newcommand{\ak}{\a_k}

\newcommand{\rk}{\R_k}
\newcommand{\ck}{\C_k}
\newcommand{\hk}{\H_k}
\newcommand{\ok}{\o_k}

\newcommand{\rz}{\R_\Z}
\newcommand{\cz}{\C_\Z}
\newcommand{\hz}{\H_\Z}
\newcommand{\oz}{\o_\Z}

\newcommand{\re}{\mathtt{Re}}


\newcommand{\dem}{\noindent \underline {\bf D\'{e}monstration :} }
\newcommand{\pr}{\noindent \underline {\bf Proof :} }
\newcommand{\indic}{\noindent \underline {\bf Indication :} }
\newcommand{\rem}{\underline {\bf Remarque :} }
\newcommand{\rek}{\noindent \underline {\bf Remark :} }
\newcommand{\fin}{\begin{flushright} \vspace{-16pt}
 $\bullet$ \end{flushright}}
\newcommand{\lpara}{\vspace{-5pt} \ \\}
\newcommand{\para}{\vspace{1pt} \ \\}
\newcommand{\Para}{\vspace{15pt} \ \\}

\newcommand{\sectionplus}[1]{\section{#1} \vspace{-5mm} \indent}
\newcommand{\subsectionplus}[1]{\subsection{#1} \vspace{-5mm} \indent}


\newcommand{\dual}{{\bf v}}
\newcommand{\com}{\mathtt{Com}}
\newcommand{\rg}{\mathtt{rg}}

\newcommand{\g}{\mathfrak g}
\newcommand{\h}{\mathfrak h}
\newcommand{\n}{\mathfrak n}
\renewcommand{\u}{\mathfrak u}


\newcommand{\spec}{\mathtt{Spec}}
\newcommand{\proj}{\mathtt{Proj}}
\newcommand{\sz}{{\spec\ \Z}}
\newcommand{\p}{{\mathbb P}}
\newcommand{\A}{{\mathbb A}}
\newcommand{\pz}{\p_\Z}
\renewcommand{\O}{{\cal O}}
\newcommand{\co}{{\cal O}}
\newcommand{\F}{{\cal F}}
\newcommand{\G}{{\cal G}}


\newcommand{\ssi}{si et seulement si }
\renewcommand{\iff}{if and only if }
\newcommand{\tr}{{}^t}
\newcommand{\trace}{\mbox{tr}}
\newcommand{\scal}[1]{\langle #1 \rangle}
\newcommand{\im}{\mathtt{Im}}


\newcommand{\suiteexacte}[3]{#1 \rightarrow #2 \rightarrow #3}
\newcommand{\surmap}{\rightarrow \hspace{-.5cm} \rightarrow}
\newcommand{\limiteinverse}{\lim_\leftarrow}
\newcommand{\liste}{\

\begin{itemize}}
\newcommand{\codim}{\mbox{codim}}
\newcommand{\point}{^{ ^\bullet} \hspace{-.7mm}}
\newcommand{\X}{\mathfrak X}

\newcommand{\res}[2]{\vspace{.15cm} 

\noindent
{\bf #1 :} {\it #2} \vspace{.15cm} 

\noindent}

\newcommand{\fonction}[5]{
\begin{array}[t]{rrcll}
#1 & : & #2 & \rightarrow & #3 \\
   &   & #4 & \mapsto     & #5
\end{array}  }

\newcommand{\fonc}[3]{
#1 : #2 \mapsto #3  }

\newcommand{\fonctionratsansnom}[4]{
\begin{array}[t]{rcl}
#1 & \rightarrow & #2 \\
#3 & \mapsto     & #4
\end{array}  }

\newcommand{\directlim}[1]{
\lim_{\stackrel{\rightarrow}{#1}}    }

\newcommand{\inverselim}[1]{
\lim_{\stackrel{\rightarrow}{#1}}    }

\newcommand{\suitecourte}[3]{
0 \rightarrow #1 \rightarrow #2 \rightarrow #3 \rightarrow 0 }

\newcommand{\matdd}[4]{
\left (
\begin{array}{cc}
{} #1 & {} #2  \\
{} #3 & {} #4
\end{array}
\right )   }

\newcommand{\matddr}[4]{
\left (
\hspace{-.2cm}
\begin{array}{cc}
{} #1 & {} \hspace{-.2cm} #2  \\
{} #3 & {} \hspace{-.2cm} #4
\end{array}
\hspace{-.2cm}
\right )   }

\newcommand{\mattt}[9]{
\left (
\begin{array}{ccc}
{}  #1 & {} #2 {} & #3 \\
{}  #4 & {} #5 {} & #6 \\
{}  #7 & {} #8 {} & #9
\end{array}
\right )   }

\title{Geometry over composition algebras : projective geometry}
\author{Pierre-Emmanuel Chaput\\
Pierre-Emmanuel.Chaput@math.univ-nantes.fr\\
Laboratoire de Mathématiques Jean Leray UMR 6629\\
2 rue de la Houssinière - BP 92208 - 44322 Nantes Cedex 3
}
\maketitle

\begin{center}
{\bf Abstract}
\end{center}

The purpose of this article is to introduce
projective geometry over composition algebras :
the equivalent of projective spaces and Grassmannians over them are
defined. It will follow from
this definition that the projective spaces are in 
correspondance with Jordan algebras and that
the points of a projective space correspond to rank one matrices in the 
Jordan algebra. A second part thus studies properties 
of rank one matrices. Finally, subvarieties of projective spaces are discussed.

\lpara

\noindent
{\it AMS mathematical classification \/}: 14N99, 14L35, 14L40. 

\begin{center}
{\bf Introduction}
\end{center}

This paper initiates a wider study of
geometry over composition algebras. The general philosophy of this study
is to discuss to what extent classical algebraic 
geometry constructions generalize
over composition algebras. Let $k$ be a commutative
field. Let $\rk,\ck,\hk,\ok$ denote the four split composition
algebras \cite{jacobson_hurwitz}.

The usual projective algebraic varieties 
over $k$ are thought of
varieties over $\rk$, and I want to 
understand analogs for $\ck,\hk$ and even
$\ok$. For example, in this article, I study ``projective spaces''
over composition algebras.

In general, in the 
octonionic case, we find varieties homogeneous under an exceptional
algebraic group. For example, the
$\ok$-generalization of the usual projective plane $\p^2_k$
is homogeneous under a group of type $E_6$ over $k$. The
analogy between this somewhat mysterious $E_6$-variety and a well-understood
projective plane allows one to understand better the geometry of this variety, 
as well as some representations of this group, which can be thought as
$SL_3(\ok)$. I plan to show that we can similarly think of $E_7$ as
``$Sp_6(\ok)$''.

\lpara

Along with this geometric and representation-theory 
motivation for studying varieties defined over
composition algebras, there is an algebraic one. Namely, the varieties that
we will meet in this context can be defined in terms of algebraic
structures (such as Jordan algebras, structurable algebras, and exceptional
Lie algebras). As I want to show, the products on these algebras correspond to
maps defined naturally in terms of geometries over a composition
algebra. This gives new insight on these algebras.

I thank Laurent Manivel for many discussions and Laurent Bonavero for comments
on a previous version of this article. I also thank Nicolas Ressayre and
Michel Brion for discussions relative to the fifth part.



\para

This paper is organised as follows : the first section recalls well-known facts
about composition algebras. A short geometric proof of the triality
principle is given. In the second section, Grassmannians over
composition algebras are defined : different definitions of $G_\a(r,n)$
as sets of $\a$-submodules of $\a^n$ with some properties are compared.

Any ``projective space'' $G_\a(1,n)$ is then seen to live in the
projectivisation of a Jordan algebra : we have naturally, as algebraic
varieties over $k$,
$G_\a(1,n) \subset \p V$, where $V$ is a Jordan algebra. Moreover,
$G_\a(1,n)$ is the variety of rank one elements in the Jordan algebra.
This is the topic of the third section : different possible definitions of
``rank one'' elements in a Jordan algebra are compared (see theorem 
\ref{th_rg1}), and the connection with the structure group of the
Jordan algebra is described.

In the fourth section, I introduce the notion of ``$\a$-subvarieties'' of
a projective space $G_\a(1,n)$. I show that there are very few of them.

Finally, section five deals with the octonionic projective plane. This
plane is defined and its automorphism group is shown to be a simple
group of type $E_6$ in all characteristics (see theorem \ref{e6}).
The projective plane over the
octonions with {\it real} coefficients $\o_\R$ has been studied extensively
\cite{tits,freudenthal}. Here, I consider the different case of the algebra
$\ok$ containing zero-divisors, and explain the new point of view of the
generalized Veronese map (theorem \ref{exceptionnel}).

\tableofcontents

\sectionplus{Background on composition algebras}

The split composition algebras have a model over $\Z$
\cite{jacobson_hurwitz} : 
the ring $\hk$ is the ring of $2 \times 2$-matrices with integral
entries. The norm of a matrix $A$ is $Q(A)=\det A$ and the conjugate of
$A=\matdd{a_{1,1}}{a_{1,2}}{a_{2,1}}{a_{2,2}}$ is the matrix
$\overline A = \matdd{a_{2,2}}{-a_{1,2}}{-a_{2,1}}{a_{1,1}}$. The
rings $\cz$ and $\rz$ are respectively the subrings of diagonal and
homothetic matrices. The ring $\oz$ may be constructed via Cayley's
process : it is the ring of couples $(A,B)$ of matrices with product
$(A,B)*(C,D)=(AC-\overline DB,B\overline C+DA)$, conjugation
$\overline{(A,B)}=(\overline A,-B)$ and norm $Q(A,B)=Q(A) + Q(B)$.
Therefore, $\rz$ and $\cz$ are commutative. For $k$ a
field and $\a \in \{ \R,\C,\H,\o \}$, we set $\ak = \az \otimes_\Z k$.
We note $\scal{x,y}=Q(x+y)-Q(x)-Q(y)$ and $\re(z)=\scal{z,1}$
(therefore $\re(1)=2$).

\begin{nota}
A composition algebra over a unital commutative ring $R$ is one of
the following algebras : $\R_R,\C_R,\H_R,\o_R$.
\end{nota}

In the sequel, many arguments will use the fact that the isotropic linear
spaces for $Q$ can be described using the algebra; namely, for $z \in \a$, 
we denote $L(z)$ (resp. $R(z)$) the image
of the left (resp. right) multiplication by $z$ in $\a$, denoted $L_z$
(resp. $R_z$).
\begin{prop}
Let $k$ be any field and
let $\a$ be a composition algebra over $k$ different from $\rk$. Let
$z,z_1,z_2 \in \a-\{0\}$ with $Q(z_1)=Q(z_2)=0$.
\begin{itemize}
\item
If $Q(z) \not = 0$, then $L(z)=R(z)=\a$. 
\item
If $Q(z)=0$, then $L(z)$ and $R(z)$ are maximal (ie of
dimension $\dim_k\a / 2$) isotropic linear subspaces of $\a$. They belong to
different connected components of the variety of maximal isotropic subspaces. 
\item
$z_2 \in L(z_1) \Longleftrightarrow z_1 \in L(z_2) 
\Longleftrightarrow \overline z_1z_2=0.$
\end{itemize}
\label{composition_general}
\end{prop}
\noindent
For example, this proposition implies that if 
$0 \not =z_1,z_2 \in \hk$
and $Q(z_1)=Q(z_2)=0$, then the dimension of $L(z_1) \cap L(z_2)$ is either
2 (ie $L(z_1)=L(z_2)$) or 0, depending on the fact that 
$z_2 \in L(z_1)$ (or equivalently
$z_1 \in L(z_2)$) or not.\\
\pr
Let $\alpha = \dim \a / 2$. The composition algebras
are alternative \cite{schafer}, which means that
$\forall x,z \in \a,z(zx)=(zz)x$. Therefore, 
$\overline z(zx)=Q(z).x$, or $L_{\overline z} \circ L_z = Q(z).Id$.

Thus, $Q(z)\not =0$ \iff $L_z$ is invertible \iff $R_z$ is. Moreover, the
kernel of $L_z$ are elements $t$ such that $R_t$ is not invertible;
therefore it is included in the quadric $\{Q=0\}$. Since
$Q(zx)=Q(z).Q(x)$, if $Q(z)=0$, then we have also $L(z) \subset \{Q=0\}$.
Since $\dim L(z)+\dim \ker L_z = \dim \a = 2\alpha$ and an
isotropic subspace has maximal dimension $\alpha$, it follows that
$\dim L(z)=\dim \ker L_z=\alpha$. The rest of the proposition is easy.
\qed

\lpara

The octonionic case is related to the triality principle. 
The following proposition
was proved in \cite{springer} using a description of $Spin_8$ involving
octonions.
\begin{prop}
If $\a = \ok$, 
the maps $L$ and $R$ induce isomorphisms from
the projective 6-dimensional quadric defined by $Q$ and the two 
connected components of the Grassmannian of
maximal isotropic spaces. Let $x,y \in \ok$ such that $Q(x)=Q(y)=0$.
\begin{itemize}
\item
$\dim (L(x) \cap L(y)) \geq 2 \Longleftrightarrow 
 \dim (R(x) \cap R(y)) \geq 2 \Longleftrightarrow \scal{x,y}=0$. 
\item
If $\dim L(x) \cap L(y) = 2$, then
$L(x) \cap L(y)=L_x[L(\overline y)]=L_y[L(\overline x)]$, and 
$R(x) \cap R(y)=R_x[R(\overline y)]=R_y[R(\overline x)]$.
\item
$xy=0 \Longleftrightarrow \dim L(x) \cap R(y) = 3$.
\end{itemize}
\label{trialite}
\end{prop}
\pr
Let $G^\pm(4,\ok)$
denote the two connected components of the variety of
maximal isotropic subspaces
of $\{Q=0\}$ in $\ok$.

First, let $x_0=\left ( \matddr1000 ,  \matddr 0000 \right )$. We have
$L(x_0)=\left ( \matddr **00 , \matddr *0*0 \right )$ and
$R(x_0)=\left ( \matddr *0*0 , \matddr 0*0* \right )$. We thus have
$L(x_0) \cap R(x_0) = k.x_0$. Since for
any $x,z \in \{Q=0\}$ with $x \not =0$ and $z \not = 0$, 
$\dim (L(x) \cap R(z)) \in \{1,3\}$, for
generic $x,z \in \{Q=0\}$, one has
$\dim (L(x) \cap R(z))=1$, so $L(x) \cap R(z) = k.xz$.

Let $x \in \{Q=0\}$ be such that for generic $z \in \{Q=0\}$,
$L(x) \cap R(z) = k.xz$. By the following lemma \ref{def_app_lin_quad}, if
$y$ is such that $L(x)=L(y)$, then there exists $\lambda \in k$ such
that $L_x=\lambda.L_y$, and so $x=\lambda.y$.

Therefore, the fiber of $[L]:\p\{Q=0\} \rightarrow G^+_Q(\ok)$ over
$L(x)$ is only $\{[x]\}$ : $[L]$ is generically injective. It follows
that it is surjective, and the same holds for $[R]$.

The previous argument is therefore valid for any $[x] \in \p\{Q=0\}$, and
$[L]$ and $[R]$ are injective. We will now show that $[L]$ is an
isomorphism. If $x,y\in \{Q=0\}$ and $y\in L(x)$, $[y] \not = [x]$, then 
$L(x) \cap L(y) = Vect(x,y)$. Therefore, given $L(x) = \Lambda$, the point 
$[x] \in \p \ok$ may be constructed as the intersection of the
projective lines 
$\Lambda \cap L(y)$, for $y \in \Lambda$. This describes the inverse
of $[L]$, which is therefore algebraic.

\lpara

I have shown that $[L]$ and $[R]$ are isomorphisms. If 
$\dim\ (L(x) \cap R(y)) = 3$, then $xy=0$, because otherwise the 
rational map
$$
\fonctionratsansnom{G^+(4,\ok) \times G^-(4,\ok)}{\p\ok}
{(\Lambda^+,\Lambda^-)}{\Lambda^+ \cap \Lambda^-}
$$ 
would be defined at $(L(x),R(y))$. Fixing $x$,
$\{y:\dim\ (L(x) \cap R(y)) = 3\}$ and $\{y:xy=0\}$ are isomorphic with
$\p^3_k$, so they are equal, proving the third point.

If $0\not =z \in L(x) \cap L(y)$, and $\lambda,\mu \in k$, then
$\overline z(\lambda x + \mu y)=0$. Therefore, 
$Q(\lambda x + \mu y)=0$ and $\scal{x,y}=0$. Since the variety 
$Q_x=\{y\in\{Q=0\}:\dim (L(x) \cap L(y)) \geq 2\}$ is a
hyperplane section of the quadric $\{Q=0\}$, we deduce
$Q_x=\{y:Q(y)=0 \mbox{ and } \scal{x,y}=0\}$.

Finally, the identity $z(\overline zt)=Q(z)t$ yields
$x (\overline yt)=\scal{x,y}t-y (\overline xt)$, which implies that
$L_x(L(\overline y)) = L_y(L(\overline x)) = L(x) \cap L(y)$, if
$\scal{x,y}=0$.
\qed

\begin{lemm}
Let $V,W$ be $k$-vector spaces, $X \subset V$ a variety and 
$f,g: V \rightarrow W$ linear maps. Assume :
\begin{itemize}
\item
$\forall x \in X- \ker f,\exists \lambda \in k : g(x)=\lambda f(x)$.
\item
No quadric of rank four vanish on $X$.
\item
$f$ has rank at least 2.
\end{itemize}
Then, $\exists \lambda \in k : g = \lambda f$.
\label{def_app_lin_quad}
\end{lemm}
\pr
Choose a basis of $W$ and let $f_1,\ldots,f_n,g_1,\ldots,g_n$ be
linear forms such that $f=(f_i),g=(g_i)$. Then the minors
$f_ig_j-f_jg_i$ vanish along $X$, so they vanish on $V$. The lemma
follows.
\qed



\sectionplus{Grassmannians over composition algebras}


\subsectionplus{Definition}

An element in the Grassmannian of $r$-planes on $k(=\rk)$ is a $k$-vector
space of dimension $r$. But if $\dim_k \a \geq 2$, $\a$ is not a field, and
this definition does not make sense anymore. However, we can define
$G_\a(r,n)$ as the set of all free right $\a$-submodules of $\a^n$ of rank $r$,
that is the set of submodules $M$ of the form:
$$
M:=\{ \sum_{t=1}^r v_t \lambda_t, \lambda_t \in \a \},
$$
with dimension $\dim \a.r$ over $k$ ($(v_t)$ is a $r$-uplet 
of elements in $\a$). 
The freeness condition generalizes the fact that the zero vector in
$k^n$ has no image in $\p^{n-1}$.

Another possible definition, which gives a structure of closed
variety, considers right $\a$-submodules of $\a^n$ of the right dimension:
$$
\tilde G_\a(r,n) = \{ E \subset \a^n : \dim_k E = \dim \a.r \mbox { and } 
\forall \lambda \in \a, E.\lambda \subset E \}.
$$

\para
We shall now study properties of these two sets. If $V$ is a 
$k$-vector space and $r$ an integer, I denote $G(r,V)$ the Grassmannian of
$r$-subspaces of $V$. Let $G(r,n)$ denote $G(r,k^n)$. I want to show
the following propositions:
\begin{prop}
$G_\H(r,n)=\tilde G_\H(r,n) \subset G(4r,\hk^n)$ is a smooth subvariety
isomorphic with the usual Grassmannian $G(2r,2n)$.
\label{g_h}
\end{prop}

\begin{prop}
The subvariety $G_\C(r,n)$ of $G(2r,\ck^n)$ is isomorphic with\linebreak
$G(r,n) \times G(r,n)$. Moreover, $\tilde G_\C(r,n)$ is the union of
$\min\{n+1-r,r+1\}$ connected components which are irreducible, one of which
is $G_\C(r,n)$.
\label{g_c}
\end{prop}
\noindent
First, let us show the following lemma (If $x$ is a real number, I
denote $[x]^+$ the least integer greater or equal to $x$):

\begin{lemm}
Let $E\subset \hk^n$ such that 
$\forall \lambda \in \hk, E.\lambda\subset E$. Then there exist
$c =  [\frac{\dim E}{4}]^+$
vectors $v_1,\ldots,v_c \in \hk^n$ such that
$E = \oplus \{v_i.\lambda,\lambda \in \hk\}$.
\label{engendre}
\end{lemm}
\rek
Assuming the result,
let $1 \leq t \leq c$ be an integer and 
$(v_{t,u})_{1\leq u\leq n}$ the coordinates of the vector $v_t$.
The kernel of the map $\hk \rightarrow \hk^n,\lambda \mapsto v_t.\lambda$ is
trivial if one of the $v_{t,u}$ is invertible, and is
$\cap_j L(\overline v_{i,j})$ otherwise, by proposition
 \ref{composition_general}. Thus, it has even dimension by
proposition \ref{composition_general}, and the rank theorem shows that
$\{ v_t. \lambda:\lambda \in \hk \}$ has also even dimension.
Thus, proposition \ref{engendre} shows that the dimension of such an $E$
is even.

\noindent
\pr
A large part of this proof holds for $\ck$; for the moment $\a$ stands
for $\ck$ or $\hk$, and accordingly we set $\alpha$ equals 1 or 2. I will
precise which argument needs $\a=\hk$.

Let $\pi_t:\a^n \rightarrow \a$ be the projection on the $t$-th
factor, and $p_t$ the restriction of $\pi_t$ to $E$. 
If $T$ is a subset of $\{1,\ldots,n\}$, let $p_T$ denote
the products of the $p_t$'s for $t \in T$.

If there exists $t$ such that $p_t$ has maximal rank $2\alpha$, then, 
choosing a vector $v$ in $E$ such that $v_t=1$, we get an isomorphism
$$
\fonction{s}{\hk \oplus (E \cap \{ x_t=0 \})}{E}{(\lambda,x)}{v.\lambda+x,}
$$
so that we are done by an inductive argument. We thus suppose that no
projection $p_t$ has maximal rank. Since $\im(p_t)$ is preserved by right
multiplication by $\a$, by proposition \ref{composition_general},
it has rank 0 or $\alpha$. In the first case, we can
also use an inductive argument, therefore, we can suppose that any projection
$p_t$ has rank $\alpha$.

This implies that for any couple $(t,u) \in \{1,\ldots n\}^2$, 
$F=p_{\{t,u\}}(E)$ has
dimension $\alpha$ or $2\alpha$. In fact, $p_{t|F}$ has rank $\alpha$ and its
kernel is preserved by right multiplication. If $t$ and $u$ are such that
$p_{\{u,t\}}$ has rank $\alpha$, then
$p_{u|p_{\{u,t\}}(E)}$ is bijective, so that if $\iota$ denotes 
$\{ 1,\cdots,n \} - \{ i\}$, $p_\iota$ is injective, and again we conclude by
an inductive argument.

It is therefore sufficient to consider the case where any projection
$p_{\{t,u\}}$ has rank $2\alpha$. Let $t$ and $u$ be arbitrary. Since
$F=p_{\{t,u\}}(E)$ is preserved by right multiplication, and since each
projection has rank $\alpha$, there exist $x,y \in \a$ such that
$F \subset L(x) \times L(y)$; since $\dim F=2\alpha$, we have equality.

The following argument works only for $\a=\hk$. If $z \in L(y)$, then
$L(z)=L(y)$, since we have $L(z) \subset L(y)$ by associativity.
Moreover, since $\a=\hk$, we can choose $z \in L(y)$ 
such that $z$ and $y$ are not
proportional; this implies $R(z) \not = R(x)$. Thus, eventually replacing
$y$ by $z$, we can assume
$R(y) \not = R(x)$. This implies by conjugation
$L(\overline y) \not = L(\overline x)$. Thus these spaces are supplementary
and the mapping
$\a \rightarrow F, \lambda \mapsto (x,y).\lambda$ is injective, proving that
$(x,y)$ is a generator of $F$. It is enough to consider a vector which
projection under $p_{\{t,u\}}$ is $(x,y)$.
\qed

\lpara
\noindent
\underline{\bf Proof of proposition \ref{g_c}:}
Let $E \subset \ck^n$ be preserved by right-multiplication.
Let us consider the base $e=\matddr 1000,f=\matddr 0001$ of $\ck$.
If $v=(v_1,\ldots,v_n)=(v_1^+e+v_1^-f,\ldots,v_n^+e+v_n^-f) \in E$,
then $v^+=(v_1^+e,\ldots,v_n^+e)=v.e \in E$, and 
$v^-=(v_1^-f,\ldots,v_n^-f)=v.f \in E$.
If $E^+$ (resp. $E^-$) denotes $E \cap (\ck^n.e)$ (resp. $E \cap (\ck^n.f)$),
then $E = E^+ \oplus E^-$.

If
$r^+$ and $r^-$ are integers between $0$ and $n$ with sum $2r$, let
$\tilde G_{r^+}(r,n)$ denote the set of linear spaces of the form
$E=E^+ \oplus E^-$, with $E^+ \subset \ck^n.e,E^- \subset \ck^n.f$ and
$\dim E^\pm = r^\pm$. Such a linear space is preserved by muliplication
by $e$ and $f$, thus it is an element of $\tilde G_\C(r,n)$. The variety 
$G_{r^+}(r,n)$ is isomorphic to
$G(r^+,n) \times G(r^-,n)$.

On the other hand, we have seen that
$\tilde G_\C(r,n) = \cup \tilde G_{r^+}(r,n)$. To prove
that the $\tilde G_{r^+}(r,n)$ are the connected components of
$\tilde G_\C(r,n)$, I recall that
for $d \in \N$, $\{ E : \dim (E \cap (\ck^n.e)) \geq d \}$
and $\{ E : \dim (E \cap (\ck^n.f)) \geq d \}$ are closed subsets of
$G(2r,\C^n_k)$.

It remains to check that $\tilde G_r(r,n)=G_\C(r,n)$. If
$(v^\pm_1,\ldots,v^\pm_r)$ is a basis of $E^{\pm}$, then $(v^+_te + v^-_tf)_t$
is a family of vectors which generates $E$.
\qed

\lpara
\underline{\bf Proof of proposition \ref{g_h}:}
Lemma \ref{engendre} says that $G_\H(r,n)=\tilde G_\H(r,n)$. In the
same way as proposition \ref{g_c}, one proves that any $E \in G_\H(r,n)$ 
may be written as $E^+ \oplus E^-$, with $E^+ \subset \hk^n.e$ and
$E^- \subset \hk^n.f$. Moreover, let $h=\matddr 0110$;
right-multiplication with $h$ is an involutive linear automorphism of
$E$ which exchanges $E^+$ and $E^-$. Therefore, they have the same 
dimension $2r$.

Moreover, it implies that giving $E$ is equivalent to
giving $E^+$, and thus the map $E \mapsto E^+$ is the desired isomorphism
between $H_\H(r,n)$ and $G(2r,2n)$.
\qed


\subsectionplus{Duality}

In this paragraph, I show an analog of the well-known fact that 
$$G(r,V) \simeq G(\dim V-r,V^*).$$

\begin{defi}
Let $V$ and $W$ be right (resp. left) $\a$-modules.
A right-linear (resp. left-linear)
map from $V$ to $W$ is a map $f:V \rightarrow W$ 
such that $\forall x,y \in V,\forall \lambda,\mu \in \a,
f(x.\lambda+y.\mu)=f(x).\lambda+f(y).\mu$
(resp. $f(\lambda.x+\mu.y)=\lambda.f(x)+\mu.f(y)$.

\noindent
A right- (resp. left-)form on $V$ is a 
right- (resp. left-)linear map from $\a^n$ to $\a$.
\end{defi}
A map $f:\a^n \rightarrow \a^m$ is
right-linear
(resp. left linear) \iff there exists a matrix $(a_{t,u})$ such that
$f((x_u))=(\sum_u a_{t,u}x_u)_t$ (resp. $f((x_u))=(\sum_u x_ta_{t,u})_t$).
Therefore, if $V$ is a free right $\a$-module of rank $n$, so is the set of
left-linear forms on $V$, which I will denote by $V^*$.

Generalizing the construction of the previous section, 
if $V$ is a free right $\a$-module
of rank $n$, let $G_\a(r,V)$ denote 
the algebraic variety parametrizing the
free right $\a$-submodules of $V$ of rank $r$. This is obviously a variety
isomorphic with $G_\a(r,n)$. Moreover, we have the following :
\begin{prop}
Let $V$ be a free right $\a$-module of rank $n$.
There is a canonical isomorphism $G_\a(r,V) \simeq G_\a(n-r,V^*)$.
\label{g(r,v*)}
\end{prop}
\pr
If $Y \subset V$ is any set, then
$Y^\bot:=\{l \in V^*: \forall y\in Y, l(y)=0\}$ is a
$k$-linear subspace of $V^*$, preserved by right multiplication
by $\a$. If $Y$ is
an element in $G_\a(r,n)$, then it is generated by $r$ vectors, and a form
vanishes on $Y$ \iff it vanishes on the generators. Therefore, $Y^\bot$ is
of dimension at least $\dim_k \a.(n-r)$. Since the pairing
$\a^n \times {\a^n_r}^*:(x,l) \mapsto \scal{1,l(x)}$ 
is perfect, the
dimension of $Y^\bot$ is exactly $\dim_k \a.(n-r)$.

If $\a=\ck$, then 
$\dim (Y^\bot \cap (\ck^n.e))=\dim (Y^\bot \cap (\ck^n.f))=n-r$.
We thus have $Y^\bot \in G(n-r,{\a^n_r}^*)$. Since
the map $L \mapsto L^\bot$ is an isomorphism, the proposition is proved.
\qed



\sectionplus{Properties of rank one matrices and \\generalized Veronese map}

This section is concerned with the particular case $r=1$ of the previous 
section, namely, I want to study projective spaces. We will see that there
is a closed connection between these spaces, Jordan algebras, and a particular
map which generalizes alltogether the Veronese map and
the quotient rational map $k^{n+1} \dasharrow \p^n$.


\subsectionplus{Background on quadratic Jordan algebras}

The satisfactory notion of Jordan algebras over a unital commutative
ring $A$ which includes the
characteristic two case is the notion of quadratic Jordan algebras.
They are by definition the free $A$-modules $V$ of finite type
with a non-zero distinguished element denoted $1$ and
equipped with a map $R:V \rightarrow End(V)$ which satisfies the following five
axioms \cite[definition 3,p.1.9]{jacobson_quadratique} :
\begin{itemize}
\item $R$ is quadratic.
\item $R(1)=Id_V$.
\item $R(a) \circ R(b) \circ R(a) = R[R(a).b]$.
\item $V_{a,b} \circ U_b = U_b \circ V_{b,a}$, if 
$V_{a,b}(x)=R(x+b).a-R(x).a-R(b).a$.
\item The two last identites hold after any extension of the base
ring $A$.
\end{itemize}

Let $\az$ be a composition algebra over $\Z$, let
$r$ be an integer, and consider the $\Z$-module $H_r(\az)$ of
hermitian matrices with entries in $\az$. If $\az$ is associative, set 
$R(A).B = ABA$, for any $A,B \in H_r(\az)$. If $\az = \oz$, then $R$
can be defined as the unique quadratic map such that
$R(A).B=ABA$ if all the coefficients of $A$ and $B$ belong to an
associative subalgebra of $\oz$. 

By \cite[theorem 5,p.1.45]{jacobson_quadratique},
if $\az$ is associative or $r \leq 3$, then the triple
$(H_r(\az),Id,R)$ is a quadratic Jordan algebra over $\Z$.
If $k$ is any field, set $H_r(\ak)=H_r(\az) \otimes_\Z k$; it is a
quadratic Jordan algebra over $k$. Refering to the ``second structure
theorem'' \cite[p.3.59]{jacobson_quadratique}, 
it is seen that these examples of
quadratic Jordan algebras play a major role in the theory of Jordan algebras.

Now, let, for
$a=1,2,4$, $V^n_a$ be the vector space defined by
$$
V^n_a = \left \{
\begin{array}{cl}
{\cal S}_n(k) & \mbox{if }a=1\\
{\cal M}_n(k) & \mbox{if }a=2\\
{\cal AS}_{2n}(k) & \mbox{if }a=4,
\end{array}
\right .
$$
where $\cal S,M,AS$ respectively stand for the set of symmetric, arbitrary,
and\\
skew-symmetric with zero diagonal entries matrices.
We choose any invertible $I\in V^n_a$, and set
\begin{equation}
R(A).B=AI^{-1}BI^{-1}A.
\label{alg_vna}
\end{equation}
It is well-known that the algebra $V^n_a$ is
isomorphic with $H_n(\a)$ (if $\dim \a=a$); 
the following corollary \ref{iso_scorza_quaternion} gives a geometric 
understanding of this isomorphism when $a=4$. In the sequel,
$I$ will be chosen to be the usual identity matrix for $a=1,2$, and the
bloc-diagonal matrix with non-vanishing coefficients $\matdd{0}{-1}{1}{0}$
for $a=4$.

We note 
$\p_\a^n:=G_\a(1,n+1)$. First, we remark that 
$\p_\C^{n-1} \simeq \p^{n-1} \times \p^{n-1} \subset \p V^n_2$
and $\p_\H^{n-1} \simeq G(2,2n) \subset \p V^n_4$ 
are naturally embedded in the projectivisations of quadratic
Jordan algebras.
The same holds for $\p_\R^{n-1}$, which embeds 
via the second Veronese embedding
in $\p V^n_1$. As we shall see in the last section, the exceptional 
quadratic Jordan
algebra $H_3(\ok)$ corresponds to the octonionic projective plane. 
The aim of this section
is to give several equivalent caracterizations of the elements in $\p V^n_a$
which correspond to points of $\p_\a^{n-1}$.


\subsectionplus{Definition of Jordan rank one matrices}

I start with a pure Jordan-algebra definition whose geometrical meaning will
become clearer in the next subsections. In $H_r(\ak)$, set
$$\scal{A,B}=\sum_{1\leq i<j\leq r} \scal{A_{i,j}B_{i,j}} + 
\sum_i A_{i,i}B_{i,i} 
\mbox{\ \ \   and\ \ \  } \trace(A) = \sum_i A_{i,i} = \scal{1,A}.$$
\begin{defi}
Let $V = H_r(\ak)$ and $0 \not = A \in V$. We will say that
the Jordan rank of $A$ is one if 
$\forall B\in V, R(A)B=\langle A , B \rangle A$.
\label{def_rang1}
\end{defi}
\rek Since $\trace$ and $\scal{.,.}$ can be defined in any quadratic Jordan
algebra, the above definition makes sense not only in $H_r(\ak)$, but in any 
(quadratic) Jordan algebra.

Each element $B \in V$ determines $\dim V$ quadratic equations given by the
coordinates of $R(A)B-\langle A , B \rangle A$. It is clear that a 
quadratic Jordan
algebra isomorphism induces an isomorphism of varieties of rank one elements.

\begin{nota}
Let ${\cal Q}_2$ denote the space of quadrics generated,
for all $B \in V$, by the coordinates of the
equation on $A$: $R(A)B-\langle A , B \rangle A=0$.
\end{nota}

Before explaining what Jordan
rank one
elements are in $V^n_a$, I want to show that this definition is well-behaved
with respect to a ``big'' algebraic group.

The structure group of a Jordan algebra $V$,
denoted $Str(V)$ is defined as the group of $g \in GL(V)$ such that
$\forall B \in V, R(g.A)=g \circ R(A) \circ \tr g$, 
the transposition being taken
with respect to the scalar product $\scal{.,.}$.
This definition may seem rather abstract to a reader not used to Jordan
algebras, so I recall that the connected component of $Str(V^n_a)$
is $GL_n,(GL_n \times GL_n)/k^*,$ or $GL(2n)$,
according to $a=1,2,4$, where $k^*$ is diagonaly embedded in 
$GL_n \times GL_n$ and the actions on $V^n_a$ are the natural ones.

\begin{lemm}
The algebraic variety of rank one elements is preserved by $Str(V)$, as well
as the vector space of quadrics ${\cal Q}_2$.
\label{rg1_str}
\end{lemm}
Therefore, a straightforward computation shows that the set of Jordan rank one
matrices
is the closed orbit of $Str(V^n_a)$ in $\p V^n_a$, namely
the set of (usual) rank one matrices if $a=1$ or 2 
and the set of rank 2 matrices
if $a=4$. Moreover, the equations ${\cal Q}_2$ are,
respectively, the two by two minors, and the Plücker equations of the 
Grassmannian. Recalling propositions \ref{g_h} and \ref{g_c},
the variety of rank one elements is thus $\p^{n-1}_\a \subset \p V^n_a$.

\noindent
\pr
Let $f \in {\cal Q}_2$; we can assume that there
exist $B,C \in V$ such that $f(A)=\scal{R(A).B-\scal{A,B}A,C}$. From the
definition of $Str(V)$, it follows that $\forall g \in Str(V)$,
$$
\begin{array}{rcl}
(g^{-1}.f)(A) & = & \scal{R(g.A).B-\scal{gA,B}gA,C}\\
         & = & \scal{g[R(A).(\tr gB)]-\scal{A,\tr gB}gA,C}\\
         & = & \scal{R(A).(\tr gB)-\scal{A,\tr gB}A,\tr gC}.
\end{array}
$$ 
Therefore this is a quadric in the vector space ${\cal Q}_2$. Thus this vector
space, and the variety it defines, are preserved by $Str(V)$.
\qed


\subsectionplus{Jordan rank and generalized Veronese maps}

In this paragraph, I give an analog of the map $k^{n+1} \dasharrow \p^n$
for $\p_\a^n$ in terms of a map which generalizes the
usual Veronese map in the case $\a=\rk$.

Following F.L. Zak \cite[th. 4.9]{zak}, let's consider the following map:
$$
\begin{array}{rccl}
\nu_2: &    {\cal A}^n & \dasharrow & \p H_n({\cal A})\\
       &    (z_i)      & \mapsto      & (z_i\overline {z_j}).
\end{array}
$$
As I will show, this map can be interpreted
as the rational map $\a^n \dasharrow \p \a^{n-1}$:
\begin{prop}
Let $\lambda \in \a$ with $Q(\lambda) \not = 0$
and $(z_i)$ such that $\nu_2(z_i)$ is well-defined.
Then $\nu_2((z_t.\lambda))$ is also well-defined and equals $\nu_2((z_t))$.
\label{nu2_projectif}
\end{prop}
\pr
In fact, $((z_t \lambda)).(\overline \lambda \overline{z_u})= 
Q(\lambda)z_t \overline{z_u}$.
\qed

\lpara
Now, let us see that the image of this map is $\p^{n-1}_\a$ :
\begin{prop}
The image of $\nu_2:\a^n \dasharrow \p H_n(\a)$ is the set of Jordan rank one
elements.
\label{veronese}
\end{prop}
\noindent
The image of this rational map is the set-theoretical one, namely the
set of all the matrices which may be written as $\nu_2(z_t)$. The proposition
shows that this set is closed.

\noindent
\pr
First, a direct computation shows
\begin{equation}
A=(a_{t,u}) \mbox{ and } B=\matdd{1}{0}{0}{0} 
\Longrightarrow ABA=(a_{t,1}a_{1,u})_{t,u}.
\label{ABA}
\end{equation}
If $A=(z_t\overline {z_u})$ is in the image of $\nu_2$,
then for $B=\matdd{1}{0}{0}{0}$, 
$ABA=Q(z_1) A$, thus is equal to
$\langle A,B \rangle A$, by (\ref{ABA}). We can also make a similar
computation for
$B=\matdd{\matdd{0}{\overline{x}}{x}{0}} {0}{0}{0}$. In fact, if
$A=(z_t \overline{z_u})$, then
$ABA=([z_t\overline{z_1}][\overline x(z_2
\overline{z_u})] + [z_t \overline{z_2}][x(z_1\overline{z_u})])_{t,u}$. Using
associativity, this equals
$$(z_t \re(\overline z_1 \overline x z_2) \overline{z_u})_{t,u}
= (z_t \re(\overline x z_2 z_1) \overline{z_u})_{t,u}
=\langle A,B \rangle (z_t\overline{z_u})_{t,u}.$$ Using the permutations and
linearity, $A$ has Jordan rank one.

In the other way, if there exists a diagonal matrix $B$ such that
$\langle A,B \rangle \not =0$, then if for example
$B=\matdd{1}{0}{0}{0}$, we see from (\ref{ABA}) that $A$ is collinear with
$\nu_2(\overline{a_{1,t}})$. If such a $B$ does not exist, then for all
diagonal $B$, $ABA=0$. Thus, using (\ref{ABA})
$$\forall t,u,v, \ \overline{a_{v,t}} a_{v,u}=0. $$ 
Then we can
assume there exists $B$ of the form
$B=\matdd{\matdd{0}{\overline x}{x}{0}} {0}{0}{0}$ such that
$\langle A,B \rangle \not = 0$; we deduce that $A$ is proportional
to $\nu_2(\overline{a_{1,t}} + \overline{a_{2,t}} x)$.
\qed

\lpara

Let $X \subset \p H_n(\hk)$ be the variety of rank one elements in the
Jordan algebra $H_n(\hk)$. Recall that $e = \matddr1000 \in \hk$.
Any $z \in \hk$ induces, by left multiplication,
a linear morphism
$R(e) \rightarrow R(e)$. Chosing the basis 
$\left \{e,\matddr 0010 \right \}$ of
$R(e)$, we can associate to $z$ a two-by-two matrix $M(z)$ representing
$L_z$. This map yields
a map $\tilde M:H_n(\hk) \rightarrow {\cal M}_{2n}$, given by
$\tilde M((a_{t,u})_{t,u})=(M(a_{t,u}))_{t,u}$. Recall that $I$
is the $2n \times 2n$ bloc-diagonal matrix with 
diagonal entries $\matdd{0}{1}{-1}{0}$.
Let $\tilde \phi$ denote the map $H_n(\hk) \rightarrow {\cal M}_{2n},
A \mapsto I.\tilde M(A)$. We also define, for $A \in H_n(\hk)$, the map
$L_A:\hk^n \rightarrow \hk^n$ defined, if $A=(a_{t,u})$, by
$L_A(z_u)=(\sum_u a_{t,u}z_u)_t$.
The following corollary gives a geometric understanding
of the isomorphism $H_n(\hk) \simeq V_4^n$.
\begin{coro}
The map
$$ \fonction{\phi}{X \subset \p (H_n(\hk))}
{G(2,k^{2n}) \subset \p (\Lambda^2 k^{2n})}{A}
{(\im L_A) \cap R(e)^n} $$
is induced by the linear isomorphism
$\tilde \phi$.
\label{iso_scorza_quaternion}
\end{coro}
\pr
The proof is a computation left to the reader (the details
are written in my thesis \cite{these}).
\qed

It is easily checked that, for the Jordan product on $V^n_4$ given by
(\ref{alg_vna}),
the map $\tilde \phi$ is also a Jordan algebra isomorphism. 
This should not come
as a surprise, as the following proposition shows.
\begin{prop}
Let $V_1$ and $V_2$ be quadratic Jordan algebras isomorphic with some algebra
$H_r(\ak)$, $I_1$ and $I_2$ their units and
$X_1 \subset \p V_1$ and $X_2 \subset \p V_2$ the corresponding 
varieties of rank one elements. Let
$f:V_1 \rightarrow V_2$ be a linear map such that $f(I_1)=I_2$. The following
conditions are equivalent:
\begin{enumerate}
\item
$f$ induces an isomorphism of varieties between $X_1$ and $X_2$.
\item
$f$ is a Jordan algebra isomorphism between $V_1$ and $V_2$.
\end{enumerate}
\end{prop}
\pr
This result will not be used in the sequel, so this proof, which uses a
little of Jordan algebra theory, could be skipped.
Of course, (2) implies (1), since the varieties $X_i$ are defined using only
the Jordan algebra structure. Now, in the quadratic Jordan algebras I consider,
there is a well-defined polynomial called the norm 
\cite{jacobson_generic_norm}. In $V^n_a$ with $a=1,2$, it
is the usual determinant of matrices and in $V^n_a$, it is the pfaffian. In the
exceptional case $H_3(\ok)$, the norm is defined by formula 
(\ref{determinant_octave}) in subsection \ref{op2}.
Moreover, I use the fact that the hypersurface defined by the norm
is the closure of the set of
sums of $r-1$ rank one elements. In fact, for the classical Jordan algebras, 
this is just an easy result of linear algebra, whereas for the exceptional
algebra, it follows from proposition \ref{exceptionnel}. Denoting $\det_1$
and $\det_2$ the norms of the Jordan algebras $V_1$ and $V_2$,
we therefore have
$\det_1(A_1)=\det_2[f(A_1)]$. Since the scalar product is the second 
logarithmic differential of the
determinant at the identity $(\scal{A,B}=D^2_I \log \det(A,B))$ 
\cite{mccrimmon}, it
follows that $\scal{A_1,B_1}_1=\scal{f(A_1),f(B_1)}_2$. Moreover, since the
quadratic product itself is also the second logarithmic derivative
of the determinant
($\scal {R(A)^{-1}.B,C}=-D^2_A \log \det (B,C)$), we deduce that $f$ is an
algebra morphism.
\qed
\para

We now relate two other possible definitions of rank one matrices to the
previous definition.
In the case $n=3$, the following proposition shows that Jordan rank one 
matrices are defined by minors (which is not the case in general).
\begin{prop}
A Hermitian matrix
$(a_{t,u})_{1 \leq t,u \leq 3}$ with value in $\a$
has Jordan rank one \iff
$$
\begin{array}{c}
a_{1,1}a_{2,2}-Q(a_{1,2})=a_{1,1}a_{3,3}-Q(a_{1,3})=a_{2,2}a_{3,3}-Q(a_{2,3})=0
\mbox{ and }\\
a_{1,1}a_{2,3}-a_{2,1}a_{1,3}=a_{3,2}a_{2,1}-a_{3,1}a_{2,2}=
a_{2,1}a_{3,3}-a_{2,3}a_{3,1}=0.
\end{array}
$$
\label{ordre_3}
\end{prop}
\pr
To prove this result, I use the corresponding one concerning the exceptional
Jordan algebras (proposition \ref{exceptionnel}, which proof is 
self-contained).
If $(a_{t,u})$ has rank one, it equals $\nu_2((z_t))$ and these minors
vanish. If these minors vanish, considering this matrix as a matrix
with coefficients in $\ok$, it has rank one by proposition \ref{exceptionnel}
which means that $\forall B \in H_3(\ok)$,
$R(A)B=\scal{A,B}A$. Thus this equality holds for
$B \in H_3(\a)$, and $A$ has rank one.
\qed

\lpara
Recall that for $A$ a matrix of order $n$ with coefficients in $\a$, 
we defined the map $L_A:\a^n \rightarrow \a^n$ by
$L_A((z_u))= (\sum_u a_{t,u}z_u)_t$. This is a $k$-linear map. If $\a$ were
a field, it would be clear that $\dim_k \a$ would divide $\dim_k \im L_A$. 
Here, it
is not the case, take for example $A=\matdd{z}{0}{0}{0}$ with $z$
a zero divisor. However, this is
true for Hermitian matrices:

\begin{prop}
Let $A \in H_n(\a)$. Then $\dim_k \a$ divides $\dim_k \im L_A$.
\label{imagediva}
\end{prop}
\pr
Suppose first that $\a=\hk$. Recall that I denote $e = \matddr 1000$ and
$f=\matddr 0001$. Under the ismoorphism of proposition 
\ref{iso_scorza_quaternion}, a matrix $A \in H_n(\hk)$ identifies with
$I.\tilde M(A)$, where $I$ stands for the bloc diagonal matrix with entries
$\matddr{0}{1}{-1}{0}$, and $\tilde M(A)$ is the
matrix of the restriction of $L_A$ to $R(e)^n$. Since $I.\tilde M(A)$ is
skew-symmetric, $\tilde M(A)$ has even rank, 
and since by associativity $R(f)^n=R(e)^n.\matddr0110$, the
rank of $L_A$ is a multiple of 4.

Considering a matrix with entries in $\ck$ as a matrix with entries in $\hk$,
we deduce the case $\a = \ck$ from the case $\a = \hk$, since
for a matrix $A$ with coefficients in $\ck$, we have
$L_A(\hk^n)=L_A(\ck^n) \oplus L_A(\ck^n).\matddr 0110$.
\qed

\begin{prop}
$A \in H_n(\a)$ has Jordan rank one \iff $L_A$ has rank $\dim_k \a$.
\label{rang_app_lin}
\end{prop}
\pr
I thank Laurent Manivel for the simplification of the argument he suggested
to me. Using the same argument as for the previous proposition, it is
enough to consider the case when $\a = \hk$.

Since for $A = \nu_2(z_t), \im L_A \subset 
\{ (z_t.\lambda),\lambda \in \a \}$ and since this rank is a multiple of
$\dim \a$, we have equality and one implication is proved. For the reverse
implication, we may by an inductive argument, left to the reader, 
suppose that $A$ has order three. If $a_{1,1} \not = 0$, the hypothesis
implies that all the columns of $A$ are right-multiple of the first, which
implies that $A$ has rank one. Therefore, we may assume that the diagonal of 
$A$ vanishes.

Moreover, if $A$ has a vanishing row, since it is Hermitian,
we in fact have to study a matrix of
the form $\matdd{0}{\overline z}{z}{0}$, with $z \in \a$ such that
$Q(z)=0$, and this matrix has rank 
one by proposition \ref{ordre_3}. 
We therefore assume that no row of $A$ vanishes.

Let $C_u$ denote the columns of $A$; I claim that all vector spaces 
$\{ C_u.\lambda,\lambda \in \a \} \subset \a^3$ 
have dimension $\dim \a/2$. In fact,
if $\{ C_u.\lambda,\lambda \in \a \}$ has dimension $\dim \a$, then the
hypothesis implies that all the columns of $A$ belong to this vector space;
since $a_{u,u}=0$, $A$ would have a vanishing row. 

Let $u$ be fixed.
Since $\{ C_u.\lambda,\lambda \in \a \}$ has dimension $\dim \a/2$, there
exists a vector space $K_u$ of dimension two such that 
$\forall t,\ker L_{a_{u,t}} \supset K_u$. Since $K_u$ is of the form
$L(z_u)$, we have $\forall t,a_{u,t} \in R(\overline z_u)$.
Since $A$ is Hermitian,
$a_{u,t}$ therefore belongs to $L(z_t) \cap R(\overline z_u)$, so it is of
the form $z_t.b_{t,u}.\overline z_u$.
The proposition \ref{ordre_3} shows that $A$ has rank one.
\qed

Finally, I would like to mention the following result, which makes a link
between my definition of rank one and another definition that we find in the
litterature \cite[p.290]{harvey}:
\begin{prop}
Let $A\in H_3(\a)$.
Then $A$ has Jordan rank one \iff $A^2=(\trace A).A$.
\end{prop}
\pr
Using definition \ref{def_rang1} with $B=Id$, we see that if $A$ has rank one,
then $A^2=(\trace A).A$ (even if we are in $H_n(\a)$ with $n>3$).

Conversly, a direct computation shows that
$\trace L_A=(\dim \a)(\trace A)$.
To prove the proposition, we can assume $\a=\hk$. 
If $\trace A=0$, $A^2=0$, so $L_A^2=0$, and so
$L_A$ has rank at most 6. By propositions \ref{imagediva} and 
\ref{rang_app_lin}, $A$ has Jordan rank one. If $\trace A=1$, $A^2=A$, so
$L_A^2=L_A$. So $L_A$ has eigenvalue 1 with multiplicity 4 and eigenvalue 0
with multiplicity 8 ($\trace L_A=4$). Therefore, the rank of $L_A$ is four
and proposition \ref{rang_app_lin} applies.
\qed


\subsectionplus{Summary of properties of rank one matrices}

Let's summarize some of the results of the preceeding subsection:
\begin{theo}
Let $V$ be a quadratic
Jordan algebra isomorphic with $H_n(\a)$ and $0 \not = A \in V$.
The following conditions are equivalent :
\begin{enumerate}
\item
$A$ has Jordan rank one.
\item
The class of $A$ in $\p V$ belongs to the closed orbit of $Str(V)$.
\item
For any isomorphism $\varphi:V \simeq {\cal S}_n(\C)$ 
(resp. ${\cal M}_n(\C)$, ${\cal AS}_{2n}(\C)$), $\varphi(A)$ has minimal
rank 1 (resp. 1,2).
\item
For any isomorphism $\varphi : V \simeq H_n(\a)$,
$\varphi(A)$ is the Veronese image of a $n$-uple of elements in $\a$.
\item
For any isomorphism $\varphi : V \simeq H_n(\a)$, $L_{\varphi(A)}$ 
has rank $\dim \a$.
\end{enumerate}
\label{th_rg1}
\end{theo}

Proposition \ref{nu2_projectif} shows that
the map $\nu_2:\a^n \dasharrow \p_\a^{n-1}$ is
exactly the analog of the map $k^n \dasharrow \p^{n-1}$. Let us understand
better this map. In the case $\a=\ck$, any vector $z \in (\ck)^n$ can be
written uniquely
as $x+y$, with $x\in R(e)^n$ and $y \in R(f)^n$. If $x$ and $y$
don't vanish, then $\nu_2(z)$ identifies via $\p H_n(\ck) \simeq \p V_2^n$ with
$([x],[y]) \in \p R(e)^n \times \p R(f)^n 
\simeq \p^{n-1} \times \p^{n-1}$. Similarly, let $h=\matddr 0110$;
a vector $z \in (\hk)^n$ can
be written uniquely as $x+y.h$, with $x,y \in R(e)^n$. The map 
$\hk^n \dasharrow G(2,R(e)^n)$ corrresponding to $\nu_2$ sends
$z=(x+y.h),x,y \in R(e)^n$ 
on the line $(x,y) \in G(2,R(e)^n)$.

Any vector $z \in \a^n$ yields a $k$-linear map 
$\a \rightarrow \a^n,\lambda \mapsto z.\lambda$. Let $\cal I$ denote the closed
subset of $\a^n$ where this map is not injective; we have:
\begin{prop}
The indeterminacy locus of the rational map 
$\nu_2:\a^n \dasharrow \p \a^{n-1}$
is exactly $\cal I$. For $z \not \in {\cal I}$,
$\nu_2^{-1}[\nu_2(z)]=\{ z.\lambda \in \a: Q(\lambda) \not = 0 \} 
\simeq \{ \lambda \in \a : Q(\lambda) \not = 0 \}.$
\end{prop}
\noindent
Therefore, we are as close as possible to the situation of a usual projective
space. In fact, if we consider the usual map $\pi : k^{n+1} \dasharrow \p^n$
(case $\a = \rk$), then $\pi(z)$ is defined \iff 
$z \not = 0$, which is equivalent to
$z \not \in {\cal I}$. Moreover, if $\pi(z)$ is
defined, then $\pi^{-1}[\pi(z)]=\{ z.\lambda : \lambda \not = 0 \}$.

\noindent
\pr
If $(z_t) \in \a^n$ is such that $\forall t,u,z_t \overline z_u = 0$, then
$\forall t,Q(z_t)=0$, and if for example $z_1 \not = 0$, then
$\forall t,z_t \in R(z_1)$. Therefore, $(z_t).\overline z_1=0$ and 
$(z_t) \in \cal I$. Thus this indeterminacy locus is included in $\cal I$.

In the case $\a=\ck$, from the preceeding description 
of the map $\nu_2$, it is clear that this
map cannot extend to a point $z$ in $R(e)$ or $R(f)$. Similarly, let
$\a=\hk$ and suppose
$z \in \hk^n$ is such that $\lambda \mapsto z.\lambda$ is not
injective. If we write as before $z=x+y.h$, 
with $x,y \in R(e)$, we deduce that $x$ and $y$ are
proportional, because $\dim \{z.\lambda,\lambda \in \hk\}=2$ implies
$\dim (\{z.\lambda,\lambda \in \hk\} \cap R(e))=1$
(proof of proposition \ref{iso_scorza_quaternion}). 
Therefore, again, we cannot extend the rational map $\nu_2$ to
$z$.
\qed 
 
As we will see in the final section,
things are not so well behaved as far as octonions are concerned.



\sectionplus{Subvarieties of projective spaces}

In this section, we assume $k$ has characteristic 0.
I define and
classify subvarieties of projective spaces: as we shall see, this 
classification is a little disappointing, because there are very few such
subvarieties. I hope that further investigations will explain this fact.

If $X \subset \p \a^n$
is a subvariety, I denote $\tilde X$ the closure of its preimage by $\nu_2$
(this preimage is a subset of $\a^{n+1}$).

From the theorem of Newlander and Nirenberg, it follows that
a {\it real} subvariety of a usual complex variety is a complex
subvariety
\iff each tangent space (which is a real subspace) is stable by the 
multiplication by complex numbers. Therefore, it seems natural to me to
introduce the following
\begin{defi}
A subvariety $X \subset \p \a^n$ is an $\a$-variety if the corresponding
affine cone $\tilde X \subset \a^{n+1}$ has the property that
$\forall x \in \tilde X,T_x\tilde X \subset \a^{n+1}$ 
is preserved by right multiplication by any element in $\a$.
\label{sous-variete}
\end{defi}
Since $\tilde X$ itself is preserved by right multiplication by $\a$, 
this condition is 
equivalent to the fact that $T_x\tilde X=T_{x.\lambda}\tilde X$.

Let us discuss what $\a$-varieties are.
\begin{prop}
Let $X \subset \p^n_\C=\p^n \times \p^n$ be a closed subvariety.
The following conditions are equivalent:
\begin{enumerate}
\item
$X$ is a $\ck$-variety.
\item
There exist $P_1,\ldots,P_a \in \ck[X_0,\ldots,X_n]$ such that
$$x \in \tilde X \Longleftrightarrow \forall 1\leq t \leq a,P_t(x)=0.$$
\item
There exist $X_1,X_2 \subset \p^n$ ordinary varieties such that
$$X = X_1 \times X_2 \subset \p^n \times \p^n.$$
\end{enumerate}
\label{ss_variete_c}
\end{prop}
\pr
I will show that (1) and (2) are equivalent to (3). 
Let us show that (3) implies (1) and
(2). If $X=X_1 \times X_2$, and if $\tilde X_i \subset k^{n+1}$ is the usual
affine cone of $X_i$, then 
$\tilde X \simeq \tilde X_1 \oplus \tilde X_2$, under the isomorphism
$(\ck)^{n+1}=R(e)^{n+1} \oplus R(f)^{n+1} 
\simeq k^{n+1} \oplus k^{n+1}$. Let us denote by $p_1$ and $p_2$ the two
projections on $k^{n+1}$, which correspond to the multiplications by $e$ 
and $f$. We have that 
$T_{(x_1,x_2)}\tilde X=T_{x_1}X_1 \oplus T_{x_2}X_2$ is preserved by $p_i$.
We thus have (1).

Let $(f_t)$ and $(g_u)$ be defining equations 
of $X_1$ and $X_2$ in $\p^n$. The inclusion of algebras $k \subset \ck$
induces an inclusion $\varphi$ of
$k[X_0,\ldots,X_n]$ in $\ck[X_0,\ldots,X_n]$. If
we set $P_t=e\varphi(f_t)$ and $Q_u=f\varphi(g_u)$, then it is clear
that $\forall (a_v,b_v) \in \C^{n+1}$, $P_t[(a_v(e)+b_v(f))_v]=0$ \iff
$f_t(a_v)=0$, and similarly $Q_u[(a_ve+b_vf)_v]=0$ \iff
$g_u(b_v)=0$. Therefore, the polynomials $(P_t,Q_u)$ define the variety
$\tilde X_1 \oplus \tilde X_2$, and (2) is true.
Note that we can reverse this computation :
if (2) is true, then $X$ is a product and we have (3).

The last thing to be proved is that (1) implies (3). If
$T_{(x_1,x_2)}\tilde X$ is stable by the two projections $p_1$ and $p_2$, 
according to the more general following proposition, $\tilde X$ is the
sum of two varieties.
\begin{prop}
Let $V_t,1\leq t\leq a$ be vector spaces over $k$ and 
$X \subset \bigoplus V_t=:V$
an irreducible affine variety such that 
$\forall (\lambda_t)\in k^a, (x_t)\in X \Rightarrow (\lambda_t.x_t) \in X$ and
$\forall x\in X, T_xX$ is the sum of its intersections with $V_t \subset V$.

Then there exist irreducible affine varieties $X_i \subset V_i$ such that 
$X=\bigoplus X_i$.
\end{prop}
\pr
By induction, we can assume that $a=2$. Let us
consider the restrictions $p_{1,2}$ to $X$ of the projections on $V_{1,2}$. 
Let $n_i$ be the
generic dimension of $T_xX \cap V_i$ and $X_i=p_i(X)$. 
Let $x=(x_1,x_2)\in X$ be a
smooth generic point. Then $T_xX=(T_xX \cap V_1) \oplus (T_xX \cap V_2)$ 
has dimension $n_1+n_2$, so $\dim X=n_1+n_2$.
Moreover, the kernel of $dp_1$ is $T_xX \cap V_2$, of dimension $n_2$, and also
the kernel of $dp_2$ has dimension $n_1$. It follows that 
$\dim X_i=n_i$. Thus $\forall x_2 \in X_2, p_2^{-1}(x_2)$ has dimension
at least $n_1$. Since the restriction of $p_1$ to this preimage is an
isomorphism on its image which is a subvariety of $X_1$, we deduce that
$p_2^{-1}(x_2)=\{ (x_1,x_2):x_1 \in X_1\}$. Therefore $X=X_1 \oplus X_2$.
\qed

\lpara

We now consider the quaternionic case:

\begin{prop}
Let $X \subset \p \hk^n=G(2,2n+2)$ be a closed subvariety.
The following conditions are equivalent:
\begin{enumerate}
\item
$X$ is a $\hk$-variety.
\item
There exist $l_1,\ldots,l_a$ right-linear forms such that
$x \in \tilde X \Longleftrightarrow \forall 1\leq t \leq a,l_t(x)=0$.
\item
There exists a linear subspace $L \subset k^{2n+2}$ such that
$X = G(2,L) \subset G(2,k^{2n+2})$.
\end{enumerate}
\end{prop}
This proposition, though a little disappointing, does not really come as a
surprise, at least from a heuristic point of view. In fact, due to the lack
of commutativity of $\hk$, the condition
that a polynomial on $\hk$ vanishes is well-defined on a quaternionic 
projective space only if this polynomial is linear.

\noindent
\pr
If $X \subset G(2,k^{2n+2})$ is the subgrassmannian $G(2,L)$, 
then $\tilde X \subset (\hk)^{n+1}$ is a sublinear
space preserved by right mutiplication by $\hk$ (in fact 
$L \subset R(e)^{n+1}$ is preserved by right 
multiplication by 1 and $e$, and
so $L \oplus L.h$ is preserved by right multiplication by $\hk$). Therefore
(3) implies (1).
By proposition \ref{g(r,v*)}, any linear space in $(\hk)^n$ preserved by
right-multiplication by $\hk$ is given by linear equations, so (3) and
(2) are equivalent.

To prove that (1) implies (3), I suggest two proves, one studying $\tilde X$ 
and one studying $X$. For the first proof, we note that 
$\tilde X \subset \hk^{n+1}$ is stable by right multiplication by 
$\ck \subset \hk$ and all
the tangent spaces $T_x\tilde X$ also. Therefore, by proposition
\ref{ss_variete_c}, there exists 
$X_1 \subset R(e)^{n+1}$ and $X_2 \subset R(f)^{n+1}$ such that 
$\tilde X=X_1 \oplus X_2$. Since $\tilde X.h=\tilde X$, $X_2=X_1.h$.
Moreover, if $x_1,x_2 \in X_1$, then $(x_1,x_2.h) \in \tilde X$ and the fact
that $T_{(x_1,x_2)}\tilde X$ is preserved by multiplication by $h$ implies
that $T_{x_2}X_2=T_{x_1}X_1.h$. Thus $T_{x_1}X_1$ does not depend on $x_1$ and
$X_1$ is a linear space.

For the second proof, we remark that if for 
$\tilde x\in \tilde X \cap (\hk^n-{\cal I})$,
$T_{\tilde x}\tilde X$ is a $\hk$-linear
space of dimension $4\dim X$ corresponding to 
$M \in G(2\dim X,R(e)^{n+1})$, 
then for $x=\nu_2(\tilde x)$, $T_{[x]}X$ is included in 
(and thus equal to by dimension count) 
$L_x^* \otimes M/L_x$, where $L_x$ stands for the linear space parametrized
by $x \in G(2,R(e)^{n+1})$. This will imply that $X$ is a
subgrassmannian by the following more general proposition.

\begin{prop}
Let $r\geq 2$ and $X \subset G(r,V)$ a subvariety of a Grassmannian such that 
$\forall x \in X$ there exists a linear subspace $M$ of $V$ such that
$T_xX=L_x^* \otimes (M/L_x) \subset L_x^* \otimes V/L_x=T_xG(r,V)$. Then there
exists a linear subspace $L \subset V$ such that $X=G(r,L)$.
\end{prop}
\pr
The hypothesis implies that $\dim X$ is a multiple of $r$, so let $d$ such
that $\dim X = rd$.
Consider the total space $T$ of the restriction of the tautological bundle to
$X$, and the natural projection $p:T \rightarrow V$. Let $L=p(T)$;
if $0 \not = v \in p(T)$, then $p^{-1}(v) \simeq \{x\in X:v \in L_x\}$.
Therefore, the tangent space to this fiber is the set of 
$\varphi \in L_x^* \otimes (V/L_x)$ such that $\varphi(v)=0$. The hypothesis
implies that this has dimension $r(d-1)$, thus $p(T)$ has dimension
$rd+r-r(d-1)=r+d$. If $v \in p(T)$, let $X_v=\{x\in X:v \in L_x\}$, and let
$T_v$ be the restriction of $T$ to $X_v$. Again, a fiber of the projection
$p:T_v \rightarrow V$ has dimension $r(d-2)$ and thus $p(T_v)$ is also of
dimension $r+d$. So, we have that $p(T_v)=p(T)$ and it follows that $p(T)$ is
a linear space, and by dimension count $X=G(r,p(T))$.
\qed



\sectionplus{The exceptional case}

In this section, I give a study of the octonionic projective plane
similar to that of the projective spaces over $\rk,\ck,\hk$. 
First, I have to understand the structure
group of the exceptional Jordan algebra. It has been known for very long
that the exceptional Lie groups of types $F_4$ and $E_6$
can be defined in terms
of Jordan algebras. Here however, I will describe the Chevalley
group of type $E_6$ (over the integers) using the incidence geometry of the 27
lines on a smooth cubic surface. This idea comes from \cite{faulkner} and 
\cite{lurie}.
First, I consider a degree
three polynomial, defined using the geometry of smooth cubic surfaces.
I show that the group of elements preserving
it is simple of type $E_6$ without using any
results on Jordan algebras (theorem \ref{e6}). Then I show that this 
polynomial is equivalent to the
determinant of the exceptional Jordan algebra (proposition \ref{det_octave}). 
Finally, using the known representation theory of
$E_6(k)$, I describe the octonionic plane
(theorem \ref{exceptionnel}).

\subsectionplus{Preliminary facts on smooth cubic surfaces}

Let $S\subset \p^3_\C$ be a smooth cubic surface and
let $\cal P$ denote the set of lines in $S$. Let $\cal L$ denote the set of
tritangent planes. If $p \in \cal P$ and $l \in \cal L$, I write $p \in l$
whenever the line lies in the plane. It is well known 
(\cite[section V.4.]{hartshorne} or \cite[section 10]{dolgachev}) that
there are 27 lines on $S$ and 45 tritangent planes, each of which
containing three lines.

Following
J.R. Faulkner \cite{faulkner}, let us call a 3-grid a couple of triples of
planes $[(l_1,l_2,l_3),(m_1,m_2,m_3)]$ such that the intersection of $l_i$ and
$m_j$ is a line in $S$ (the incidence relation of the 9 corresponding lines 
looks like a ``$3 \times 3$-grid'').

\subsectionplus{Definition of the Chevalley group of type $E_6$}

Following
J.R. Faulkner, let $\theta:{\cal L} \rightarrow \{-1,1\}$ be a function with
the property that for any 3-grid $[(l_1,l_2,l_3),(m_1,m_2,m_3)]$ of
$({\cal P},{\cal L})$, one has
$$\theta(l_1)\theta(l_2)\theta(l_3) + \theta(m_1)\theta(m_2)\theta(m_3) = 0$$
(theorem 5 in \cite{faulkner} exhibits such a function).

Let $V=\Z^{\cal P}$ and let $\alpha$ be the following form on this module:
$$\alpha(f)=\sum_{l \in {\cal L}} \theta(l) \prod_{p \in l} f(p).$$
We have the following:

\begin{theo}
The group-scheme of elements preserving $\alpha$ is 
isomorphic with the simply-connected 
Chevalley group of type $E_6$. Its projectivisation is the adjoint group.
Moreover, if $k$ is 
algebraically closed, then the closed orbit of $G(k)$ 
acting on $\p (V \otimes k)$
is the singular locus of the cubic hypersurface defined by $\alpha$.
\label{e6}
\end{theo}
\pr
Let $G$ be this group-scheme.
Let $k$ be an infinite field; let us first show that $G(k)$ is split
reductive of type $E_6$. 

First of all, there is an explicit formula for $\alpha$, computed as
formula (7) in \cite{faulkner}: let $V_1$ be the $\Z$-module of
$3 \times 3$-matrices with integer coefficients and 
$W=V_1 \oplus V_1 \oplus V_1$. Let $\beta$ be the cubic form
$$\beta(A,B,C)=\det(A)+\det(B)+\det(C)-\trace(ABC).$$
Then $(V,\alpha)$ is isomorphic with $(W,\beta)$.

I can exhibit a maximal torus in $G$: 
let $M,N,P\in SL_3(k)$. Then the action
$(M,N,P).(A,B,C)=(MAN^{-1},NBP^{-1},PCM^{-1})$ defines an element in $G(k)$.
Taking diagonal matrices, we thus have a torus 
$T \subset G$ of rank 6 (defined over $\Z$). Let us show that it is a
maximal torus in $G$. Let $g \in G(k)$ be an element commuting with $T(k)$.
Since $g$ preserves the eigenlines of $T(k)$, it is of the form
$(g.f)(p)=\lambda(g)f(p)$. Therefore, it preserves the three spaces
$V_1 \oplus \{0\} \oplus \{0\},\{0\} \oplus V_1 \oplus \{0\}$ and
$\{0\} \oplus \{0\} \oplus V_1$. Since $\alpha(A,0,0)=\det(A)$, 
there exist $M_1,N_1$ such that $g.(A,0,0)=(M_1AN_1^{-1},0,0)$. Moreover,
since $g$ acts diagonaly, $M_1$ and $N_1$ are diagonal. Similarly, we prove
that $g.(A,B,C)=(M_1AN_1^{-1},N_2BP_1^{-1},P_2CM_2^{-1})$. The fact that $g$
preserves $\trace(ABC)$ implies $M_1=M_2,N_1=N_2$ and $P_1=P_2$ and so
$g \in T(k)$.

I now show that the Weyl group of $G(k)$ is the Weyl group of type $E_6$.
Let $g$ be in the normalizer of $T(k)$. Then $g$ permutes the eigenlines of
$T(k)$; therefore it induces a bijection of $\cal P$. Since $g$ preserves
$\alpha$, this bijection corresponds to a bijection of the incidence
$({\cal P},{\cal L})$.
Since the isomorphism
group of this geometry is $W(E_6)$ \cite[th. 23.9]{manin},
we therefore have a
map $W(G) \rightarrow W(E_6)$. We have already seen that this map is injective.
Let us argue for its surjectivity.

Let $w$ be an isomorphism of $({\cal P},{\cal L})$.\\ Consider the function
$\fonction{\psi}{\cal L}{\{-1,1\}}{l}{\theta(l)/\theta(w.l)}$. It
has the property that 
$\psi(l_1)\psi(l_2)\psi(l_3) = \psi(m_1)\psi(m_2)\psi(m_3)$
for any 3-grid $[(l_1,l_2,l_3),(m_1,m_2,m_3)]$. By 
lemma 4 in \cite{faulkner}, there exists $x \in \{-1,1\}^{\cal P}$ such that
$\psi(l)=\prod_{p\in l} x(p)$. Therefore, we can set
$(g.f)(p)=x(p)g(w.p)$ to get an element $g \in N_G(T)$ which is
equivalent to $w$ modulo $T$. 

From this it follows that $V \otimes k$ is an
irreducible representation of $G(k)$. In fact, let 
$U \subset V \otimes k$ be any sub-representation. If $U$ contains a
vector different from $0$, since $k$ is infinite, $U$ contains an
eigenvector for $T \otimes k$. Using the above action of the Weyl
group of $G(k)$, $U$ contains all the eigenvectors and therefore
equals $V \otimes k$. Thus, $G(k)$ is reductive.

Let $G(k)^0$ be the connected component of the identity element in $G(k)$.
I have to show that the image of $N_{{G(k)}^0}(T)$ is also $W(E_6)$.
Since this
image is a normal subgroup of the group $W(E_6)$ which has a normal
simple subgroup of index 2, it is enough to
exhibit an odd element of this image. 
This is easy, considering
the action $(A,B,C) \mapsto (MAN^{-1},NBP^{-1},PCM^{-1})$ with 
$(M,N,P)$ in the connected variety $SL_3(k)^3$. Note that considering
these elements of $G(k)$, one checks that $G(k)$ and $G(k)^0$ have the
same center, namely $\{j.Id:j^3=1\}$.

\para

I therefore have shown that $G(k)^0$ is a split reductive
group of type $E_6$. To show that $G(k)$ is in fact connected, I first
prove the result about the singular locus.

Assume that $k$ is algebraically closed. Let
$X' \subset \p (V \otimes k)$ denote the closed $G(k)^0$-orbit 
and let $X$ be the singular
locus of the hypersurface defined by $\alpha$. 
Since $X$ is a non-empty closed invariant subvariety of
$\p (V \otimes k)$, we have $X \supset X'$.
Let $K$ and $K'$ be the spaces of quadrics
which vanish along $X$ and $X'$: $K \subset K'$. We have a
$G(k)^0$-equivariant map $\varphi:(V\otimes k) \rightarrow K$ given by
$v \mapsto (u \mapsto D_u\alpha (v))$. It is easy to check that $\varphi$
is not identically 0. Since $V$ is irreducible, 
by Schur's lemma, this is an injection and $\dim \varphi(V \otimes k)=27$.
Moreover, since $X' \subset \p (V \otimes k)$ is projectively normal 
\cite[theorem 1]{ramanan}, the restriction map yields an exact sequence
$$\suitecourte{K'}{Q(V\otimes k)}{H^0(X',{\cal O}(2))}.$$
$(Q(V\otimes k)$ denotes the space of quadratic forms on $V \otimes k$).
The dimensions of the vector spaces involved 
in this exact sequence are the same in positive characteristic as in
zero characteristic (the middle dimension is obviously the same and the
others could only eventually be larger, since $X'$ can be realized as a
flat scheme over $\Z$); therefore, $\dim K' = 27$. We thus have 
$\varphi(V) \subset K \subset K'$ and these vector spaces have dimension
27, so $\varphi(V)=K=K'$. Since the ideal of $X'$ is generated by quadrics 
\cite[th 3.8, p.86]{ramanathan}, $X=X'$.

\para

Let $P$ the stabilizor in $G(k)^0$ of a point in $X$. We have $X=G(k)^0/P$. 
It follows from \cite[Th\'eor\`eme 1]{demazure} that 
$G(k)^0 \rightarrow Aut(X)$ is surjective.
Since $X$ is the singular locus of $\{\alpha=0\}$, we have exact
sequences
($C(k)$ denotes the common center of $G(k)$ and $G(k)^0$)
$$
\begin{array}{c}
1 \rightarrow C(k) \rightarrow G(k) \rightarrow Aut(X) \rightarrow 1\\
1 \rightarrow C(k) \rightarrow G(k)^0 \rightarrow Aut(X) \rightarrow 1.
\end{array}
$$
Thus we have: $G(k)=G(k)^0$.

The center of $G(k)$ contains 
$\{j.Id:j^3=1\}$ which in any characteristic is a scheme of length three.
Therefore $G(k)$ is simply-connected.

Hence for all algebraically closed fields $k$, $G(k)$ is the 
simply-connected simple group of type $E_6$.
Moreover, it is proved in \cite[Theorem 5.5.1]{lurie} that the dimension
of the Lie algebra of $G$ over $\Z/p\Z$ is allways 78 (it does not depend on
$p$). It follows therefore from the proof of 
\cite[proposition III 10.4]{hartshorne}
that $G$ is smooth over $\Z$. From the uniqueness result 
\cite[expos\'e XXIII, corollaire 5.4]{sga}, the theorem follows.
\qed

\subsectionplus{$E_6$ and the exceptional Jordan algebra}

\label{op2}

The following result makes the link between the preceeding subsection and
the rest of the article. Recall that the determinant of the
exceptional Jordan algebra is defined by the equation
\cite[(18), p.37]{jacobson_generic_norm} :
\begin{equation}
\det \mattt{r_1}{\overline x_3}{\overline x_2} {x_3}{r_2}{x_1} {x_2}
{\overline{x_1}}{r_3} = r_1r_2r_3 + 2 \langle x_1x_2,x_3 \rangle
-r_1Q(x_1)-r_2Q(x_2)-r_3Q(x_3).
\label{determinant_octave}
\end{equation}

\begin{prop}
The previous form $\alpha$ is isomorphic with the determinant of the
exceptional Jordan algebra $\det$.
\label{det_octave}
\end{prop}
\noindent
Therefore, the group of elements preserving $\det$ is also the simple
simply-connected group of type $E_6$.\\
\pr
If $a_{i,j},b_{i,j},c_{i,j}, \ 1\leq i,j \leq 3,$ are integers,
a courageous reader will check that the determinant of the hermitian
matrix
$
H=\mattt{b_{1,3}}{\overline {x_3}}{\overline{x_2}} 
{x_3}{c_{3,1}}{x_1}  {x_2}{\overline{x_1}}{-a_{1,1}}
$
with
$$
\begin{array}{l}
x_1 = \matdd {a_{2,1}} {-c_{3,3}} {a_{3,1}} {c_{3,2}}
+ \matdd {b_{3,1}} {-b_{2,1}} {-b_{3,2}} {b_{2,2}} e \ \ ;\\ \\
x_2 = \matdd {a_{1,2}} {a_{1,3}} {-b_{3,3}} {b_{2,3}}
+ \matdd {c_{2,2}} {-c_{2,3}} {-c_{1,2}} {-c_{1,3}} e \ \  ;\\ \\
x_3 = \matdd {-a_{3,3}} {a_{3,2}} {a_{2,3}} {-a_{2,2}}
+ \matdd {b_{1,2}} {c_{1,1}} {b_{1,1}} {c_{2,1}} e\\
\end{array}
$$
is $\det A + \det B + \det C -\trace(ABC)$, if
$A=(a_{i,j}),B=(b_{i,j})$
and $C=(c_{i,j})$.
\qed

\para

Here is the explanation how I found this formula. A Sch\"afli's double-six
is by definition \cite[p.403]{hartshorne}
a couple $(E_i,F_i)$ of sextuples of lines in $S$
such that $E_i$ don't meet $E_j$ if $i\not =j$, 
and $E_i$ meets $F_j$ \iff $i \not =j$.
Such double-sixes exist and label uniquely the other lines in $S$,
since there is a unique line meeting $E_i$ and $F_j$ when $i\not = j$.

Now, let us say that a linear form $l \in \{a_{i,j},b_{i,j},c_{i,j}\}$
``meets'' another linear form $m$ if $lm$ divides a monomial appearing in
the expression of $\alpha$. It is easily seen that
$((a_{1,1},a_{2,1},a_{3,1},b_{2,1},b_{2,2},b_{2,3}),
(a_{1,2},a_{2,2},a_{3,2},b_{1,1},b_{1,2},b_{1,3}))$ is a Sch\"afli's
double-six for this incidence relation. Playing the same game with
$\det$, one sees that we can
start filling a hermitian matrix of coordinates with the following forms :
$$
\mattt{b_{1,3}}
{\overline{x_3}}
{\overline{x_2}}
{x_3=\matddr 0 {a_{3,2}} 0 {a_{2,2}} + \matddr {b_{1,2}} 0 {b_{1,1}} 0 e}
{c_{3,1}}
{x_1=\matddr {a_{2,1}} 0 {a_{3,1}} 0 + \matddr 0 {b_{2,1}} 0 {b_{2,2}} e}
{x_2=\matddr {a_{1,2}} 0 0 {b_{2,3}}}
{\overline{x_1}}
{a_{1,1}}
.
$$
Then, using the fact that a Sch\"afli's double-six labels all the
lines, one can finish filling the above matrix. One gets the matrix $H$
up to signs; the determinant of this matrix involves the 45 expected
monomials, but with wrong signs. These signs may be corrected 
using the algorithm
described in the proof of \cite[lemma 4]{faulkner}.

\subsectionplus{The octonionic projective space}

Let ${\cal L}$ be the set of all the projective lines
in the space
$\p \{\re = 0 \} \subset \p (\ok)$, on which the restriction of the octonionic
product vanishes identically. Let $X_0$ be the set of matrices of the form
$\mattt{0}{\overline a}{\overline b} {a}{0}{\overline c} {b}{c}{0}$, with
$a,b,c$ octonions which generate a line in $\cal L$. Let $X_1$ be the 
set-theoretic image by $\nu_2$ of triples of elements in $\ok$
generating an associative subalgebra.

\begin{prop}
Let $X\subset \p H_3(\ok)$ be the variety of rank one elements in the
exceptional Jordan algebra. Then
\begin{itemize}
\item
$X = X_0 \amalg X_1$.
\item
$X$ is the closed orbit of $Str(H_3(\ok))$.
\item
$X$ is the singular locus of the hypersurface $\{ \det =0 \}$.
\item
$X$ is defined by the following quadrics:
$$
\begin{array}{r}
a_{1,1}a_{2,2}=Q(a_{1,2}),\ a_{1,1}a_{3,3}=Q(a_{1,3}),\ 
a_{2,2}a_{3,3}=Q(a_{2,3})\\
a_{1,1}a_{2,3}=a_{2,1}a_{1,3},\ a_{3,2}a_{2,1}=a_{3,1}a_{2,2},\ 
a_{2,1}a_{3,3}=a_{2,3}a_{3,1}.
\end{array}
$$
\item
The hypersurface defined by the determinant is the closure of the set of 
sums of two rank one elements.
\end{itemize}
\label{exceptionnel}
\end{prop}
\rek
This variety is not, according to \cite[th.4.9, p.90]{zak}, the image of all
octonionic triples. In fact, the coefficients of a matrix in $X$ belong to
an associative subalgebra of $\ok$, which is not the case in general if we
take the image of any triple. Moreover, all the elements of $X$ are not
images of $\nu_2$. The claim of F.L. Zak, that 
$\nu_2(z.\lambda)=\nu_2(z)$ for any invertible octonion $\lambda$, is also
wrong, due to the lack of associativity.\\
J. Roberts \cite{roberts} has shown that the 
singular locus of the hypersurface defined by $\det$ is a Severi variety (cf 
\cite{zak,severi}
for the definition and study of Severi varieties).\\
If $A \in \p H_3(\o_\C)$ is in fact 
in $\p H_3(\o_\R)$ and annihilates the quadrics
of the proposition, its diagonal cannot vanish; therefore it is easy to see
that it is the image by $\nu_2$
of a vector in $\o_\R^3$ with one coordinate equal to 1. We thus see the link
with the set of matrices considered by Freudenthal or Tits 
\cite{freudenthal,tits}.

\noindent
\pr
We already know that the closed orbit is the singular locus of 
$\{\det = 0\}$.
A direct computation using the explicit
formula (\ref{determinant_octave}) shows
that the equations of this locus are
\begin{eqnarray}
a_{1,1}a_{2,2}=Q(a_{1,2}),\ a_{1,1}a_{3,3}=Q(a_{1,3}),\ 
a_{2,2}a_{3,3}=Q(a_{2,3}) \nonumber\\
a_{1,1}a_{2,3}=a_{2,1}a_{1,3},\ a_{3,2}a_{2,1}=a_{3,1}a_{2,2},\ 
a_{2,1}a_{3,3}=a_{2,3}a_{3,1}.
\label{quadrique_octave}
\end{eqnarray}
One can check that an element of rank one annihilates these quadrics;
therefore $X$ is the closed orbit.

Computing the number of roots of the parabolic subgroup
stabilizing a highest weight
vector, it is easily seen that $\dim X=16$.

If $X_1^0$ is the image by $\nu_2$ of vectors of the form $(1,z_1,z_2)$, it is
clear that the preceeding quadrics vanish on $X_1^0$, therefore
$\overline X_1^0 \subset X$. Since $\dim X_1^0=16$, 
equality holds. Since $\ok$
is an alternative algebra (meaning that every subalgebra generated by two
elements is associative), all coefficients of a matrix in $X'$ belong to an
associative algebra. This explains why we consider images by $\nu_2$
of triples of octonions generating an associative algebra.

If $A$ and $B$ belong to $X$, then $\det$ vanishes at order two
on the line $(AB)$ at the points $A$ and $B$. Therefore, all the points of the
line $(AB)$ are with vanishing determinant. Now, if 
$A=\mattt{1}{0}{0} {0}{0}{0} {0}{0}{0}$ and
$B=\mattt{0}{0}{0} {0}{1}{0} {0}{0}{0}$, then 
$T_AX \cap T_BX = \left \{ \mattt{0}{*}{0} {*}{0}{0} {0}{0}{0} \right \}.$
Therefore, the closure of the set of sums of two elements
of $X$ is at least an hypersurface in $H_3(\ok)$. So, this is exactly the
hypersurface defined by the determinant.
\lpara

To finish the proof of the proposition, we have to understand the image $X_1$
of $\nu_2$. Let $A \in X$. If $A$ has non-vanishing first diagonal coefficient,
we can suppose that this coefficient equals one, and then
$A=\nu_2(1,a_{2,1},a_{3,1})$ belongs to $X_1^0$ and thus $X_1$. The same holds
for any matrix which diagonal does not vanish. If $\re(a_{2,1}) \not = 0$, we
can argue as in proposition \ref{veronese} with the matrix
$\mattt{0}{1}{0} {1}{0}{0} {0}{0}{0}$. It is therefore sufficient to consider
matrices of the form
$\mattt{0}{\overline a}{\overline b} {a}{0}{\overline c} {b}{c}{0}$ with
$\re(a)=\re(b)=\re(c)=0$. In this case, the quadrics (\ref{quadrique_octave})
show that $Q(a)=Q(b)=Q(c)=0$ and $ab=ba=ac=ca=bc=cb=0$ (since 
$\overline a=-a,\overline b=-b$ and $\overline c=-c$). 
Thus the octonions $a,b,c$ are in a subalgebra
where the product vanishes identically. This subalgebra has dimension at most 2
(in fact if $a$ and $b$ are not collinear, then since
$\overline a c= \overline b c=0$, $c \in L(a) \cap L(b)=Vect(a,b)$). If it has
dimension 2, then $A \in X_0$. If it has dimension 1, then all coefficients
of $A$ are scalar multiple of some octonion with vanishing norm; since we can
put this octonion in a subalgebra of $\ok$ isomorphic with $\hk$, proposition
\ref{veronese} shows that $A \in X_1$.

\lpara

The last thing to prove is that $X_0$ and $X_1$ are disjoint.

\noindent
Let
$A=\mattt{0}{\overline a}{\overline b} {a}{0}{\overline c} {b}{c}{0} \in X_0$,
with $a,b,c$ octonions as before, and suppose there exist $z_1,z_2,z_3$ such 
that $\nu_2(z_1,z_2,z_3)=A$, or:
$$
\left \{
\begin{array}{c}
Q(z_1)=Q(z_2)=Q(z_3)=0.\\
z_1 \overline z_2=a\\
z_1 \overline z_3=b\\
z_2 \overline z_3=c.
\end{array}
\right .
$$

Suppose first that no element in $\{a,b,c\}$ is zero.

\noindent
If $a$ and $b$, $b$ and $c$, and $c$ and $a$ are not collinear, we deduce from
the system that $z_1 \in L(a) \cap L(b)=(a,b)$, and similarly
$z_2,z_3 \in (a,b)$. We therefore have a contradiction, because this implies
$z_1 \overline z_2=0$.

\noindent
If for example $a$ and $b$ are collinear, but not $a$ and $c$, then, the
system implies $z_2,z_3 \in L(a) \cap L(c)=(a,c)$. We therefore have
$z_2 \overline z_3 = 0$, a contradiction.

If for example $c=0$ and $a,b \not = 0$, we still have $z_1 \in (a,b)$, but 
only $z_2 \in L(a)$ and $z_3 \in L(b)$. If we had $z_2 \in L(b)$, we would
have $z_2 \in (a,b)$ and $a=z_1 \overline z_2=0$, a contradiction.
Thus, $z_2 \not \in L(b)$ and so $bz_2 \not = 0$. This shows
that $L(b) \cap R(z_2)=\scal {bz_2}$ (proposition \ref{trialite}). 
Since $z_3 \in L(b) \cap R(z_2)$, it is
propotionnal to $bz_2$. But in that case, since $z_2 \in L(a)$,
$\overline z_3$ is in $R_b[R(a)]=R(a) \cap R(b)=(a,b)$ by proposition
\ref{trialite}. 
This in turn would imply that $z_1 \overline z_3=0$, and a contradiction.
\qed


\begin{thebibliography}{99}
  \bibitem[Cha 02]{severi}      P.E. Chaput, {\it Severi varieties.}
Math. Z. {\bf 240} (2002), no. 2, 451--459.

  \bibitem[Cha 03]{these}        P.E. Chaput, thesis, available at\\
http://www.math.sciences.univ-nantes.fr/$\tilde{\hspace{.2cm}}$chaput

  \bibitem[DV 04]{dolgachev}  I.V. Dolgachev, A Verra, 
{\it Topics in classical algebraic geometry}, available at
http://www.math.lsa.umich.edu/$\tilde{\hspace{.2cm}}$idolga/

  \bibitem[Fre 54]{freudenthal} H. Freudenthal,
{\it Beziehungen der $E\sb{7}$ und $E\sb{8}$ zur Oktavenebene},
Indagat. Math. {\bf 16} (1954) 218--230 and 636--638, 
{\bf 17} (1955) 151--157 and 277--285, {\bf 25} (1963) 457--487.

  \bibitem[SGA 3]{sga} Grothendick, A.,Demazure, M. 
     Sch\'emas en groupes. III: Structure des sch\'emas en groupes r\'eductifs. 
     S\'eminaire de G\'eom\'etrie Alg\'ebrique du Bois Marie 1962/64 (SGA 3). 
     Lecture Notes in Mathematics {\bf 153}.
     
  \bibitem[Dem 77]{demazure} Demazure, M. {\it Automorphismes et d\'eformations des vari\'et\'es de Borel.}  Invent. Math. {\bf 39} (1977), no. 2, 179--186.
  
  \bibitem[Fau 01]{faulkner} Faulkner, J.R. {\it Generalized quadrangles and cubic forms.}  Comm. Algebra  {\bf 29}  (2001),  no. 10, 4641--4653.

  \bibitem[Gar 01]{garibaldi} Garibaldi, R.S.
{\it Structurable algebras and groups of type $E\sb 6$ and $E\sb 7$.} 
J. Algebra {\bf 236} (2001), no. 2, 651--691.

  \bibitem[Har 90]{harvey}      F.R. Harvey,
Spinors and calibrations.
Perspectives in Mathematics, 9.
Academic Press, Inc., Boston, MA, 1990.

  \bibitem[Har 77]{hartshorne} Hartshorne, R. Algebraic geometry. Graduate Texts in Mathematics, {\bf 52}. Springer-Verlag, 1977.

  \bibitem[Jac 58]{jacobson_hurwitz} N. Jacobson,
{\it Composition algebras and their automorphisms.}
Rend. Circ. Mat. Palermo (2) {\bf 7} (1958), 55--80.

  \bibitem[Jac 59]{jacobson_groupe} N. Jacobson,
{\it Some groups of transformations defined by Jordan algebras. I.} 
J. Reine Angew. Math. {\bf 201} (1959) 178--195.

  \bibitem[Jac 63]{jacobson_generic_norm} N. Jacobson,
{\it Generic norm of an algebra.}  Osaka Math. J. {\bf 15} (1963) 25--50.

  \bibitem[Jac 68]{jacobson}    N. Jacobson,
Structure and representations of Jordan algebras.
Am. Math. Soc. Coll. Publ. vol XXXIX (1968).

  \bibitem[Jac 69]{jacobson_quadratique}   N. Jacobson,
Lectures on quadratic Jordan algebras. Tata Institute of Fundamental 
Research Lectures on Mathematics, No. {\bf 45}. 
Tata Institute of Fundamental Research, Bombay, 1969

  \bibitem[Jac 85]{jacobson2}  N. Jacobson,
{\it Some projective varieties defined by Jordan algebras.}
J. Algebra {\bf 97} (1985), no. 2, 565--598.

  \bibitem[LM 99]{magique}     J.M. Landsberg, L. Manivel,
{\it The projective geometry of Freudenthal's magic square.}
J. Algebra {\bf 239} (2001), no. 2, 477--512.

  \bibitem[Lur 01]{lurie} Lurie, J. {\it On simply laced Lie algebras and their minuscule representations.}
Comment. Math. Helv. {\bf 76} (2001), no. 3, 515--575.

  \bibitem[Man 74]{manin}  Manin, Y.I. Cubic forms. Algebra, geometry, arithmetic.
    Second edition. North-Holland Mathematical Library, {\bf 4}.

  \bibitem[McC 65]{mccrimmon}   K. McCrimmon,
{\it Norms and non-commutative Jordan algebras.}
Pacific J. Math {\bf 15}(1965) 925-956.


  \bibitem[RR 85]{ramanan}    S. Ramanan, A. Ramanathan,
{\it Projective normality of flag varieties and Schubert varieties.}
Invent. Math. {\bf 79} (1985), no. 2, 217--224.

  \bibitem[Ram 87]{ramanathan}  A. Ramanathan, 
{\it Equations defining Schubert varieties and Frobenius 
splitting of diagonals.} 
Inst. Hautes \'Etudes Sci. Publ. Math. {\bf 65}, (1987), 61--90.

  \bibitem[Rob 88]{roberts}     J. Roberts,
Projective embeddings of algebraic varieties.
Monografías del Instituto de Matemáticas, 19.
Universidad Nacional Autónoma de México, México, 1988.


  \bibitem[Sch 66]{schafer}     R.D. Schafer,
An introduction to nonassociative algebras. Corrected reprint of the 1966 original.
Dover Publications, Inc., New York, 1995


  \bibitem[BS 60]{springer}    F. van der Blij, T.A. Springer,
{\it Octaves and triality.}
Nieuw Arch. Wisk. (3) {\bf 8} 1960 158--169.

  \bibitem[Tit 53]{tits}        J. Tits, 
{\it Le plan projectif des octaves et les groupes de Lie exceptionnels.}
Acad. Roy. Belgique. Bull. Cl. Sci. {\bf 39}(5), (1953). 309--329.


  \bibitem[Zak 93]{zak}         F. L. Zak,
Tangents and Secants of Algebraic Varieties.
American Mathematical Society, Providence, RI (1993).

\end{thebibliography}
\end{document}